\documentclass[a4paper,12pt]{article}
\begin{document}

\leftline{\bf \Large Did Hypatia Know about Negative Numbers? }
\vskip 0.2truecm
\noindent Hypatia lived in Alexandria in the 4th-5th century AD. She was one of the most
remarkable women in history. Her father Theon, presumably the last member of the Museum,
instructed her in mathematics, philosophy and classical Greek literature. With time she
herself became a brilliant teacher of Neoplatonic philosophy and mathematics. Numerous young
aspirants from wealthy Christian and pagan families came to Alexandria with the sole purpose
of joining the sophisticated inner circle of her students.
\vskip0.2truecm \noindent
Theon and Hyptia edited several of the major mathematical treatises available in their time:
{\it Arithmetica} by Diophantus, {\it Almagest} by Ptolemy, {\it Conics} by Appolonius of
Perga, and {\it Elements} by Euclid. A significant part of this work, in particular on the
most difficult and demanding {\it Arithmetica}, she did alone. To some of these treatises
she added commentaries and exemplary exercises.
\vskip0.2truecm \noindent
A charismatic figure in the city's intellectual life and a close friend to many important
officials, Hypatia had a significant influence on the Alexandrian elite and politics. She
was admired for her knowledge, wisdom, righteousness and personal charm. Her public lectures
for the Alexandrian intelligentsia attracted a great deal of interest. However, her elevated
social position and popularity did not prevent her from becoming the victim of a brutal
murder in 415 AD by the fanatic mob who accused her of conducting forbidden magical
practices. Groundless rumors about Hypatia's dealings with witchcraft were most likely
inspired by her political enemies. Bishop Cyril was certainly one of them, but the
allegations concerning his direct involvement have never been proven. The tragic
circumstances of her death created the romantic legend of Hypatia that persists in Western
literature, art and drama.
\vskip0.2truecm \noindent
There have been two absorbing books published recently which deal specifically with the
figure of Hypatia. Each of them, in its own way, sets to demythologize her life, legacy or
legend. In {\it Hypatia of Alexandria}\footnote{Maria Dzielska (tr. F. Lyra), {\it Hypatia
of Alexandria}. Cambridge (Mass.): Harvard University Press, 1995. ({\it Revealing
Antiquity}, 8). Pp. xi + 157. ISBN 0-674-43775-6.}, Maria Dzielska proves that the real
Hypatia was anything but a young and liberated rebel fighting against Christianity, whereas
Michael A. B. Deakin in {\it Hypatia of Alexandria: Mathematician and
Martyr}\footnote{Michael A. B. Deakin, {\it Hypatia of Alexandria: Mathematician and
Martyr}, Amherst (N.Y.): Prometheus Books, 2007; also: Michael A.B. Deakin ``Hypatia and Her
Mathematics'', {\it The American Mathematical Monthly}, March 1994, Volume 101, Number 3,
pp. 234-243.} argues that she was not a mathematical genius, and her contributions to
mathematical knowledge were slight or non-existent.
\vskip0.2truecm \noindent
In this Letter we comment on one particular aspect of Hypatia's enigmatic biography by
translating into English a short poem that appeared in a recent review of the third revised
Polish edition of Maria Dzielska's book\footnote{Marek Abramowicz, ``Niezwyk{\l}a uroda
r{\'o}wna{\'n} Diofantosa'', {{\it \'Swiat Nauki}}, October 2010, pp. 88-92.}. It poses a
simple and specific question: did Hypatia know about the negative numbers?
\vskip0.2truecm \noindent
Diophantus noticed that equations like $4 = 4X + 20$ have no solution and called them
``absurd'' or ``false''. Did Hypatia, his careful editor and commentator, share this
opinion? The Chinese and Indians had known about the negative numbers long before Hypatia,
and used them in their pure mathematics as well as in practical book-keeping: {\it a debt
cut off from nothingness becomes a credit; a credit cut off from nothingness becomes a
debt}. Knowing certainly about credit and debt, could Hypatia realize that $X = -4$ is the
solution to the ``absurd'' equation of Diophantus? Neither Dzielska, nor Deakin, nor any
other scholar, knows the answer to this question. There are no manuscripts and we have no
evidence. But the absence of evidence should not be mistaken for the evidence of absence.
\vskip0.2truecm \noindent
%
\vskip 0.2truecm
\leftline{\bf Ageometretos}
\vskip 0.2truecm
\leftline{By the shore of the uninviting sea,}
\leftline{Far away from Greece,}
\leftline{I discussed with another Hypatia}
\leftline{The Negative Numbers.}
\vskip 0.2truecm
\leftline{Now her name was Paola,}
\leftline{Her lashes superbly long,}
\leftline{And her slim hands}
\leftline{Cradled the memory of Tergeste.}
\vskip 0.2truecm
\leftline{We were amused by the petty thing,}
\leftline{A Diophantine problem,}
\leftline{Perhaps known to Hypatia:}
\leftline{Among whole numbers}
\leftline{Find all the threes}
\leftline{Of the roots of an equation}
\leftline{$X + Y + Z = 3$}
\leftline{$X^3 + Y^3 + Z^3 = 3$}
\vskip 0.2truecm
\leftline{One such three,}
\leftline{$(X, Y, Z) = (+1, +1, +1)$}
\leftline{Is apparent to all.}
\leftline{Everyone understands also}
\leftline{That in other solutions,}
\leftline{If they exist,}
\leftline{At least one of the three}
\leftline{Is a Negative Number.}
\leftline{For example}
\leftline{$(X, Y, Z) = (+4, +4, -5)$}
\vskip 0.2truecm
\leftline{Did divine Hypatia}
\leftline{Know the Negative Numbers?}
\leftline{Could she discover them}
\leftline{In the vague reference}
\leftline{Of Diophantus?}
\vskip 0.2truecm
\leftline{None of us can know it.}
\vskip 0.2truecm
\leftline{$(X, Y, Z) = (+1, +1, +1)$}
\leftline{$(X, Y, Z) = (+4, +4, -5)$}
\vskip 0.2truecm
\leftline{Are there any other}
\leftline{Matching threes?}
\leftline{How to find them All?}
\leftline{I know, and so does Paola,}
\leftline{And perhaps}
\leftline{Hypatia knew it.}
\vskip 0.2truecm
\leftline{But this cannot be deciphered}
\leftline{By those missing the Platonic gift}
\leftline{For mathematics.}
\vskip 0.2truecm
\noindent Following the idea outlined in the poem, let us assume
(without loss of generality) that $Z < 0$ and write,
\begin{equation}
\label{Z-first}
X + Y = 3 - Z
\end{equation}
\begin{equation}
\label{Z-second}
X^3 + Y^3 = 3 - Z^3
\end{equation}
After dividing Eq. (\ref{Z-second}) side by side by Eq.
(\ref{Z-first}) we arrive at,
\begin{equation}
\label{master-24}
X^2 + Y^2 -XY -Z^2 - 3Z - 9 = \frac{24}{Z - 3}.
\end{equation}
By substituting $Y = 3 - Z - X$ into Eq. (\ref{master-24}), we may
write,
\begin{equation}
\label{quadratic}
X^2 + (Z - 3)X - 3Z = \frac{8}{Z - 3}.
\end{equation}
From the above equation we see that $Z - 3$ must divide $8$. This
condition reduces the set of possible solutions for $Z$ to just
two,
\begin{equation}
\label{set-second}
Z = \left(-5, -1 \right),
\end{equation}
of which $Z =-1$ should be rejected. Thus, all possible solutions
to Hypatia's problem are,
\begin{eqnarray}
\label{all-solutions}
(X,Y,Z)~~ = &&(+1,+1,+1) \nonumber \\
&&(+4,+4,-5) \nonumber \\
&&(+4,-5,+4) \nonumber \\
&&(-5,+4,+4).
\end{eqnarray}
The above problem was presented to the competitors of the Polish 1963 {\it
Ma\-the\-ma\-ti\-cal Olympiad}, of which one of the authors of this Letter was a Finalist.
We do not know the original author of the problem.
\vskip 0.2truecm
\leftline{Marek Abramowicz}
\leftline{Physics Department}
\leftline{G{\"o}teborg University, Sweden}
\leftline{\tt marek.abramowicz@physics.gu.se}
\vskip 0.2truecm
\leftline{Anna Cetera}
\leftline{English Literature Department}
\leftline{Warsaw University, Poland}
\leftline{\tt a.cetera@uw.edu.pl}
\end{document}